\documentclass[11pt]{amsart}  
\usepackage{amsthm, amsmath, amscd, amssymb, latexsym, stmaryrd, color}

\usepackage[top=3.0cm, bottom=3.0cm, left=3.0cm, right=3.0cm]{geometry}

\theoremstyle{plain}

\newtheorem{theorem}{Theorem}[section]

\newtheorem{prop-def}[theorem]{Proposition-Definition}

\newtheorem{conjecture}[theorem]{Conjecture}

\theoremstyle{definition}

\theoremstyle{remark}

\usepackage[mathscr]{eucal}
\usepackage{graphics, graphpap}
\usepackage{array, tabularx, longtable}
\usepackage{color}
\numberwithin{equation}{section}

\def\Var{\mathrm{Var}}
\def\loc{\mathrm{loc}}
\def\sr{\mathrm{sr}}
\def\ord{\mathrm{ord}}
\def\Spec{\mathrm{Spec}}
\def\Field{\mathrm{Field}}
\def\DP{\mathrm{DP}}
\def\P{\mathrm{P}}

\def\DDef{\mathrm{Def}}
\def\RDef{\mathrm{RDef}}

\def\ac{\mathrm{ac}}
\def\L{\mathbb{L}}
\def\I{\mathrm{I}}

\title[The integral identity conjecture]{A proof of the integral identity conjecture, II}  
\author{L\^e Quy Thuong}
\date{}
\address{Department of Mathematics, Vietnam National University \newline
\indent 334 Nguyen Trai Street, Hanoi, Vietnam}
\email{leqthuong@gmail.com}

\address{BCAM - Basque Center for Applied Mathematics \newline \indent Alameda de Mazarredo 14, E-48009 Bilbao, Basque Country, Spain}
\email{qle@bcamath.org}

\thanks{This research is funded by the Vietnam National University, Hanoi (VNU) under project number QG.16.06. It is also supported by ERCEA Consolidator Grant 615655 - NMST and by the Basque Government through the BERC 2014-2017 program and by Spanish Ministry of Economy and Competitiveness MINECO: BCAM Severo Ochoa excellence accreditation SEV-2013-0323.}

\keywords{motivic integration, motivic zeta function, motivic nearby cycles, integral identity conjecture, definable subassignment, measurable subassignment}
\subjclass[2010]{Primary 03C10, 14E18, 14G10}

\begin{document}           

\begin{abstract}
In this note, using Cluckers-Loeser's theory of motivic integration, we prove the integral identity conjecture with framework a localized Grothendieck ring of varieties over an arbitrary base field of characteristic zero.
\end{abstract}

\maketitle  

\section{Statement of conjecture and main theorem}
\subsection{Equivariant Grothendieck ring of varieties}
Let $k$ be a field of characteristic zero, $S$ an algebraic $k$-variety, and $\Var_S$ the category of $S$-varieties. Let $K_0(\Var_S)$ be the Grothendieck ring of $\Var_S$, which is the quotient of the free abelian group generated by the $S$-isomorphism classes $[X\to S]$ in $\Var_S$ such that $[X\to S]=[Y\to S]+[X\setminus Y\to S]$ for any Zariski closed subvariety $Y$ of $X$. It is a commutative ring with respect to fiber product.

We consider the projective system of $\mu_n=\Spec k[t]/(t^n-1)$ with transitions $\mu_{mn}\to \mu_n$ given by $\lambda\mapsto \lambda^m$, and define $\hat\mu=\varprojlim \mu_n$. A good $\mu_n$-action on an $S$-variety $X$ is a group action each of whose orbits is contained in an affine $k$-subvariety of $X$, a good $\hat{\mu}$-action on $S$-variety $X$ is a good $\mu_n$-action for some $n$. The $\hat\mu$-equivariant Grothendieck group $K_0^{\hat{\mu}}(\Var_S)$ of $S$-varieties endowed with good $\hat\mu$-action is the quotient of the free abelian group generated by the $\hat\mu$-equivariant isomorphism classes $[X\to S,\sigma]$, $\sigma$ being a good $\hat{\mu}$-action on $S$-variety $X$, modulo the conditions $[X\to S,\sigma]=[Y\to S,\sigma|_Y]+[X\setminus Y\to S,\sigma|_{X\setminus Y}]$, for $Y$ $\sigma$-stable Zariski closed in $X$, and $[X\times_k\mathbb A_k^n\to S,\sigma]=[X\times_k\mathbb A_k^n\to S,\sigma']$, whenever $\sigma$ and $\sigma'$ lift the same $\hat{\mu}$-action on $X\to S$ to an affine action on $X\times\mathbb A_k^n\to S$. The structure of a commutative ring with unity on $K_0^{\hat{\mu}}(\Var_S)$ is given by fiber product. Let $\L$ be the class of the trivial line bundle over $S$, with trivial $\hat\mu$-action. We write $\mathscr M_S^{\hat{\mu}}$ for the localization of $K_0^{\hat{\mu}}(\Var_S)$ inverting $\L$, and $\mathscr M_{S,\loc}^{\hat\mu}$ for the localization of $\mathscr M_S^{\hat\mu}$ inverting the elements $1-\L^{-n}$, for every $n$ in $\mathbb N_{>0}$. The ring $\mathscr M_{\Spec k,\loc}^{\hat\mu}$ is rewritten simply by $\mathscr M_{\loc}^{\hat\mu}$.

Any morphism of $k$-varieties $g: S\to S'$ induces a ring morphism $g^*: \mathscr M_{S'}^{\hat{\mu}}\to \mathscr M_S^{\hat{\mu}}$ by fiber product, and induces a group morphism $g_!: \mathscr M_S^{\hat{\mu}}\to \mathscr M_{S'}^{\hat{\mu}}$ by composition. When $S'$ is $\Spec k$ we replace the symbol $g_!$ by the symbol $\int_S$. Let $\loc$ denote the natural morphism $\mathscr M_k^{\hat\mu}\to \mathscr M_{\loc}^{\hat\mu}$.
 
Consider the ring $\mathscr M_S^{\hat\mu}[[T]]$, and its subset $\mathscr M_S^{\hat\mu}[[T]]_{\sr}$ of rational series, which consists of $\mathscr M_S^{\hat\mu}$-polynomials in variables $\frac{\L^pT^q}{(1-\L^pT^q)}$, with $(p,q)$ in $\mathbb{Z}\times\mathbb{N}_{>0}$. There exists by \cite{DL1} a unique $\mathscr M_S^{\hat\mu}$-linear morphism $\lim_{T\to\infty}: \mathscr M_S^{\hat\mu}[[T]]_{\sr}\to \mathscr M_S^{\hat\mu}$ such that $\lim_{T\to\infty}\frac{\L^pT^q}{(1-\L^pT^q)}=-1$, for every $(p,q)$ in $\mathbb{Z}\times\mathbb{N}_{>0}$.

\subsection{Motivic nearby cycles}
Let $X$ be a smooth algebraic $k$-variety of pure dimension $d$. For $n\geq 1$, let $\mathscr L_n(X)$ be the $k$-scheme of $n$-jets on $X$, which represents the functor sending a $k$-algebra $A$ to the set of morphisms of $k$-schemes $\Spec(A[t]/(t^{n+1}))\to X$. These schemes together with morphisms $\mathscr L_m(X) \to \mathscr L_n(X)$ ($m\geq n$) induced by truncation form a projective system, and we denote its limit by $\mathscr L(X)$. Let $f$ be a regular function on $X$ with nonempty zero locus $X_0$. For $n\geq 1$, let $\mathscr{X}_n(f)$ be the $k$-variety of $n$-jets $\varphi$ in $\mathscr{L}_n(X)$ with $f(\varphi)= t^n\mod t^{n+1}$, which admits an obvious morphism to $X_0$ and the natural $\mu_n$-action $(\lambda,\varphi(t))\mapsto \varphi(\lambda t)$. Write $[\mathscr{X}_n(f)]$ for the class of $\mathscr{X}_n(f)\to X_0$ in $\mathscr M_{X_0}^{\hat{\mu}}$. By \cite{DL1}, the series $Z_f(T):=\sum_{n\geq 1}[\mathscr{X}_n(f)]\L^{-nd}T^n$ in $\mathscr M_{X_0}^{\hat{\mu}}[[T]]$ is rational, and the limit $\mathscr S_f:=-\lim_{T\to\infty}Z_f(T)$ in $\mathscr M_{X_0}^{\hat{\mu}}[[T]]$ is called the motivic nearby cycles of $f$. If $x$ is a closed point in $X_0$, we may consider the motivic Milnor fiber of $f$ at $x$, $\mathscr S_{f,x}=i_x^*\mathscr S_f$ in $\mathscr M_k^{\hat{\mu}}[[T]]$, where $i_x$ is the inclusion of $\{x\}$ in $X_0$.

\subsection{Conjecture and main theorem}
The integral identity conjecture plays a crucial role in Kontsevich-Soibelman's theory of motivic Donaldson-Thomas invariants for noncommutative Calabi-Yau threefolds \cite{KS}. We now state the version for regular functions of the conjecture (for the full version, see \cite[Conjecture 4.4]{KS}).

\begin{conjecture}[\cite{KS}]\label{conj}
Let $(x,y,z)$ be coordinates of the $k$-vector space $k^d=k^{d_1}\times k^{d_2}\times k^{d_3}$. Let $f$ be in $k[x,y,z]$ such that $f(0,0,0)=0$ and $f(\lambda x, \lambda^{-1} y,z)=f(x,y,z)$ for $\lambda$ in $\mathbb G_{m,k}$. Then the identity $\int_{\mathbb A_k^{d_1}}i^*\mathscr{S}_f=\L^{d_1}\mathscr{S}_{\tilde{f},0}$ holds in $\mathscr M_k^{\hat\mu}$, with $\tilde{f}$ the restriction of $f$ to $\mathbb A_k^{d_3}$, and $i$ the inclusion of $\mathbb A_k^{d_1}$ in $f^{-1}(0)$.
\end{conjecture}

Here, we identify $\mathbb A_k^{d_1}$ with $\mathbb A_k^{d_1}\times \{0\}\times \{0\}$, hence by the homogeneity of $f$, we consider it as a subvariety of $f^{-1}(0)$. We also identify $\mathbb A_k^{d_3}$ with $\{0\}\times \{0\}\times \mathbb A_k^{d_3}$, thus by definition, $\tilde{f}(z)=f(0,0,z)$.

The conjecture was first proved in \cite{Thuong1} in the case where $f$ is either a function of Steenbrink type or the composition of a pair of regular functions with a polynomial in two variables. In \cite[Theorem 1.2]{Thuong}, we show that, if the field $k$ is algebraically closed, Conjecture \ref{conj} holds in $\mathscr M_{\loc}^{\hat\mu}$. Recently, under the weaker assumption that the base field $k$ contains all roots of unity, Nicaise and Payne \cite{NP} prove the conjecture with the full context $\mathscr M_k^{\hat\mu}$.

The main result of this note is the following theorem.

\begin{theorem}\label{ThmA}
Conjecture \ref{conj} is true in $\mathscr M_{\loc}^{\hat\mu}$, namely, $\loc\left(\int_{\mathbb A_k^{d_1}}i^*\mathscr{S}_f\right)=\loc\left(\L^{d_1}\mathscr{S}_{\tilde{f},0}\right)$.
\end{theorem}

Note that our proof for the theorem does not use the assumption that $k$ is algebraically closed, that is, $k$ may be any field of characteristic zero. The materials for the proof are in Cluckers-Loeser's motivic integration of constructible motivic functions \cite{CL}.

\section{Measurable subassignments}\label{SS4}

\subsection{Definable subassignments}
We consider the formalism of Cluckers and Loeser \cite{CL} with a concrete Denef-Pas language $\mathcal L_{\DP,\P}$ consisting of the ring language $\{+,-,\cdot,0,1\}$ for valued fields, also the ring language for residue fields, and the Presburger language $\{+,-,0,1,\leq\}\cup \{\equiv_n\mid n\in \mathbb N_{>0}\}$ for value groups, where $\equiv_n$ is the equivalence relation modulo $n$. Let $\Field_k$ be the category of algebraically closed fields $K$ containing $k$ in the $\mathcal L_{\DP,\P}$-language, where sentences take coefficients in $k$ and $k(\!(t)\!)$, and morphisms of $\Field_k$ are field morphisms. The theory corresponding to $\Field_k$ is the theory of algebraically closed fields containing $k$, each model of this theory is a triple $(K(\!(t)\!),K,\mathbb Z)$ with $K$ in $\Field_k$. The valued fields $K(\!(t)\!)$ are endowed with a natural valuation map $\ord_t: K(\!(t)\!)^{\times}\to \mathbb Z$ augmented by $\ord_t(0)=+\infty$, and with a natural angular component map $\overline{\ac}: K(\!(t)\!)\to K$, with convention $\overline{\ac}(0)=0$. 

A basic affine definable subassignment has the form $h[m,n,r]$, where $h[m,n,r](K)=K(\!(t)\!)^m\times K^n\times\mathbb Z^r$. More generally, if $W=\mathcal X\times X\times\mathbb Z^r$ with $\mathcal X$ a $k(\!(t)\!)$-variety and $X$ a $k$-variety, we define $h_W(K):=\mathcal X(K(\!(t)\!))\times X(K)\times\mathbb Z^r$. An arbitrary definable subsassignment is a set of points in $h[m,n,r]$, or in $h_{W}$, satisfying a given formula $\varphi$.  

Among a broad collection of definable subassignments, we now only consider the category $\DDef_k$ of affine definable subassignments where objects are pairs $(\mathsf Z,h[m,n,r])$ with $\mathsf Z$ being a definable subassignment of $h[m,n,r]$, and a morphism $(\mathsf Z,h[m,n,r])\to (\mathsf Z',h[m',n',r'])$ is a definable morphism $\mathsf Z\to \mathsf Z'$. Due to \cite{CL}, by a definable morphism $\mathsf Z\to \mathsf Z'$ one means a morphism of subassignments $\mathsf Z\to \mathsf Z'$ such that its graph is a definable subassignment of $h[m+m',n+n',r+r']$. Let $\RDef_k$ be the full subcategory of $\DDef_k$ whose objects are definable subassignments of $h_{\mathbb A_k^n}$ for $n$ in $\mathbb N$.

Let $X$ be an affine algebraic $k$-variety. A {\it (good) $\mu_n$-action} on $h_X$ is a definable morphism of definable subassignments $h_{\mu_n\times X}\to h_X$ such that the corresponding morphism of $k$-varieties
$\mu_n\times_k X\to X$ is a (good) $\mu_n$-action. A {\it good $\hat\mu$-action} on $h_X$ is a good $\mu_n$-action on $h_X$ for some integer $n\geq 1$. For an algebraic $k(\!(t)\!)$-variety $\mathcal X$, the definable subassignment $h_{\mathcal X}$ admits a natural $\mu_n$-action $h_{\mu_n}\times h_{\mathcal X} \to h_{\mathcal X}$ induced by $(\lambda,t)\mapsto \lambda t$, for all $n\in\mathbb N_{>0}$. The profinite group scheme $\hat\mu$ acts naturally on $h_{\mathcal X}$ via $\mu_n$ for some integer $n\geq 1$. 

The Grothendieck semiring and ring of the category $\RDef_k$ are defined in \cite[Section 5.1.2]{CL}, however, in this note we only want to work with its $\hat\mu$-equivariant version. By definition, the {\it $\hat\mu$-equivariant Grothendieck group $K_0^{\hat\mu}(\RDef_k)$} is the quotient of the free abelian group generated by definable $\hat\mu$-equivariant isomorphism classes $[\mathsf X,\sigma]$, with $\mathsf X$ in $\RDef_k$ endowed with a good $\hat\mu$-action $\sigma$, modulo the relations: $[\mathsf X,\sigma]=[\mathsf Y,\sigma|_{\mathsf Y}]+[\mathsf X\setminus \mathsf Y,\sigma|_{\mathsf X\setminus \mathsf Y}]$, for $\mathsf Y$ $\sigma$-stable definable subassignment of $\mathsf X$, and $[\mathsf X\times h_{\mathbb A_k^m},\sigma]=[\mathsf X\times h_{\mathbb A_k^m},\sigma']$, whenever $\sigma$ and $\sigma'$ lift the same $\hat\mu$-action on $\mathsf X$ to an affine action on $\mathsf X\times h_{\mathbb A_k^m}$, for any integer $m\geq 0$. The cartesian product of subassignments induces a commutative with unity ring structure on $K_0^{\hat\mu}(\RDef_k)$.

Put $\mathbb A:=\mathbb Z\left[\L,\L^{-1},\frac{1}{1-\L^{-n}}\mid n\in \mathbb N_{>0}\right]$ 
where, by abuse of notation, $\L$ also stands for the class of $h_{\mathbb A_k^1}$ in $K_0^{\hat\mu}(\RDef_k)$. By quantifier elimination for the theory of algebraically closed fields containing $k$ (i.e., the Chevalley constructibility), definable subassignments of $h_{\mathbb A_k^n}$, for $n$ in $\mathbb N$, are defined by formulas without quantifiers, thus objects in $\RDef_k$ may be viewed as constructible sets. This correspondence is compatible with the $\hat\mu$-actions mentioned above. Hence, there are canonical isomorphisms of rings $K_0^{\hat\mu}(\RDef_k)\cong K_0^{\hat\mu}(\Var_k)$ and $K_0^{\hat\mu}(\RDef_k)\otimes_{\mathbb Z[\L]}\mathbb A \cong \mathscr M_{\loc}^{\hat\mu}$.

\subsection{Motivic measure}
In view of Theorem 10.1.1 in the paper \cite{CL}, there is a unique functor from $\DDef_k$ to the category of abelian semigroups, $\mathsf X \mapsto \I C_+(\mathsf X)$, which assigns to the projection $\mathsf X \to h_{\Spec k}$ a morphism of semigroups $\mu: \I C_+(\mathsf X)\to \I C_+(h_{\Spec k})$, such that the eight axioms (A1) to (A8) in that theorem, characterizing an integration theory, are satisfied. By \cite[Proposition 12.2.2]{CL}, if $\mathsf X$ is a definable subassignment of $h[m,n,0]$ which is bounded, i.e., there exists an $s$ in $\mathbb N$ such that $\mathsf X$ is contained in the subassignment of $h[m,n,0]$ defined by $\ord_t x_i\geq -s$ for $1\leq i\leq m$, then the characteristic function $\mathbf 1_{\mathsf X}$ is in $\I C_+(\mathsf X)$. In this case we call $\mathsf X$ {\it motivically measurable} and its motivic measure $\mu(\mathsf X):=\mu(\mathbf 1_{\mathsf X})\in \I C_+(h_{\Spec k})$. Also by \cite{CL}, there is a canonical morphism from $\I C_+(h_{\Spec k})$ to $K_0(\Var_k)\otimes_{\mathbb Z[\L]}\mathbb A$, thus by composition we can consider $\mu(\mathsf X)$ as an element of $K_0(\Var_k)\otimes_{\mathbb Z[\L]}\mathbb A$.

In the previous definition of boundedness, if we can take $s=0$, $\mathsf X$ is called {\it small} (see \cite[Section 16.3]{CL} for a more general definition of {\it small definable subassignments}). There is a canonical action of $\hat\mu$ on $h[m,0,0]$ induced by $(\lambda,t)\mapsto \lambda t$. We say that the definable subassignment $\mathsf X$ is stable under this action if there exists a natural number $n\geq 1$ such that, for every $x=(x_1(t),\dots,x_m(t))$ in $\mathsf X$ and $\lambda$ in $\mu_n$, the point $\lambda \cdot x=(x_1(\lambda t),\dots,x_m(\lambda t))$ is in $\mathsf X$. Since formulas defining $\mathsf X$ are in the language $\mathcal L_{\DP,\P}$, by quantifier elimination for algebraically closed fields, they also define a semi-algebraic subset $X$ of the arc space $\mathscr L(\mathbb A_k^m)$ of $\mathbb A_k^m$. The assignment $\mathsf X \mapsto X$ carries the canonical $\hat\mu$-action on $h[m,0,0]$ to the canonical $\hat\mu$-action on $\mathscr L(\mathbb A_k^m)$, and in that way, $X$ is also stable for the action on $\mathscr L(\mathbb A_k^m)$. As in \cite[Theorem 16.3.1, Remark 16.3.2]{CL} we can see that $X$ is measurable as $\mathsf X$ is measurable, and that since (with the above action) $\mu'(X)$ is in $\mathscr M_{\loc}^{\hat\mu}$, the measure $\mu(\mathsf X)$ of $\mathsf X$ is also in $K_0^{\hat\mu}(\RDef_k)\otimes_{\mathbb Z[\L]}\mathbb A \cong \mathscr M_{\loc}^{\hat\mu}$. Here, as explained in \cite[Theorem 16.3.1]{CL}, $\mu'$ stands for Denef-Loeser's motivic measure \cite{DL2}, and further by \cite[Remark 16.3.2]{CL}, we can consider that this measure takes value in $\mathscr M_{\loc}^{\hat\mu}$.

\section{Sketch of proof of Theorem \ref{ThmA}}
Let us consider the motivic zeta function $Z_f(T)$ of the polynomial $f$ in the theorem. Write the $n$-th coefficient of $\int_{\mathbb A_k^{d_1}}i^*Z_f(T)$ as follows: 
$$\int_{\mathbb A_k^{d_1}}i^*[\mathscr{X}_n(f)]\L^{-nd}=\left[\mathscr{U}_n\right]\L^{-nd}+\left[\mathscr{W}_n\right]\L^{-nd},$$ 
where $\mathscr{U}_n$ is the set of jets in $\mathscr{X}_n(f)\cap \big(\mathscr L_n(\{0\}\times \mathbb A_k^{d_2+d_3})\cup \mathscr L_n(\mathbb A_k^{d_1}\times\{0\}\times \mathbb A_k^{d_3})\big)$ originated in $\mathbb A_k^{d_1}$, $\mathscr{W}_n$ is the set of jets in $\mathscr{X}_n(f)$ that are originated in $\mathbb A_k^{d_1}$ and not contained in $\mathscr{U}_n$. The elements $\mathscr{U}_n$ and $\mathscr{W}_n$ are $\mu_n$-stable, and they give rise to rational series with coefficients in $\mathscr M_k^{\hat\mu}$. Because of the hypothesis on $f$, taking $\lim_{T\to \infty}$ for the decomposition of the $\int_{\mathbb A_k^{d_1}}i^*Z_f(T)$ reduces the proof of Theorem \ref{ThmA} to checking that $\loc\big(\lim_{T\to \infty}\sum_{n\geq 1}\left[\mathscr{W}_n\right]\L^{-nd}T^n\big)$ vanishes in $\mathscr M_{\loc}^{\hat\mu}$. The latter has been a challenging problem, and the previous attempts \cite{Thuong1}, \cite{Thuong} and \cite{NP} for solving it had to use certain additional assumptions.

Now we write $\mathscr{W}_{n,m}$ for the set of $(\varphi_1,\varphi_2,\varphi_3)$ in $\mathscr{W}_n$ with $\ord_t\varphi_1+\ord_t\varphi_2=m$, and let us observe that it is still stable under the canonical $\mu_n$-action. In what follows we only consider the set $\mathscr{W}_{n,m}$ when it is nonempty. Suggested ideally from \cite{Thuong}, with the hypothesis of Theorem \ref{ThmA}, our approach is to construct a constructible set $\widetilde{\mathscr{W}}_{n,m}$ endowed with a good $\hat\mu$-action and a $\hat\mu$-equivariant constructible surjective morphism $\mathscr{W}_{n,m}\to \widetilde{\mathscr{W}}_{n,m}$ such that its fiber over a point of residue field $k'$ is isomorphic to $\mathbb A_{k'}^{n+1}\setminus \mathbb A_{k'}^{n+1-m}$. Once we can do this, it follows (not obvious) that $[\mathscr{W}_{n,m}]=[\widetilde{\mathscr{W}}_{n,m}]\L^{n+1}(1-\L^{-m})$ in $\mathscr M_{\loc}^{\hat\mu}$, and then, Theorem \ref{ThmA} will be proved completely, because the rest of the proof is elementary. 

Our idea is that we use Cluckers-Loeser's motivic integration \cite{CL}, together to a slight development to $\hat\mu$-action context, as seen in Section \ref{SS2}. Clearly, $\mathbb G_{m,k(\!(t)\!)}$ is an algebraic group and the action of $\mathbb G_{m,k(\!(t)\!)}$ on the $k(\!(t)\!)$-variety 
$$\mathcal X:=\big(\mathbb A_{k(\!(t)\!)}^{d_1}\setminus\{0\}\big)\times_{k(\!(t)\!)} \big(\mathbb A_{k(\!(t)\!)}^{d_2}\setminus\{0\}\big)\times_{k(\!(t)\!)} \mathbb A_{k(\!(t)\!)}^{d_3}$$ 
given by 
$$\tau\cdot(\varphi_1,\varphi_2,\varphi_3):=(\tau\varphi_1,\tau^{-1}\varphi_2,\varphi_3)$$ 
is free. It follows that the space of its orbits is an algebraic $k(\!(t)\!)$-variety and the canonical projection $\phi: \mathcal X \to \mathcal Y:= \mathcal X/\mathbb G_{m,k(\!(t)\!)}$ is a surjective morphism of algebraic $k(\!(t)\!)$-varieties. This morphism $\phi$ induces a definable morphism $h_{\phi}: h_{\mathcal X} \to h_{\mathcal Y}$ of definable subassignments in the theory of Cluckers and Loeser. Take the preimage $\mathscr{W}_{n,m}^{\infty}$ of $\mathscr{W}_{n,m}$ under the canonical morphism $\mathscr L(\mathbb A_k^d)\to \mathscr L_n(\mathbb A_k^d)$. By \cite[Section 16.3]{CL}, $\mathscr{W}_{n,m}^{\infty}$ corresponds to a small definable subassignment $\mathsf V_{n,m}$ of the basic definable subassignment $h[d,0,0]$, both have the same measure $[\mathscr{W}_{n,m}]\L^{-nd}$ in $\mathscr M_{\loc}^{\hat\mu}\cong K_0^{\hat\mu}(\RDef_k)\otimes_{\mathbb Z[\L]}\mathbb A$ in Denef-Loeser's motivic measure and Cluckers-Loeser's motivic measure, endowed with $\hat\mu$-action, respectively. By Denef-Pas quantifier elimination theorem, $h_{\phi}(\mathsf V_{n,m})$ is a definable subassignment and the restriction $h_{\phi}|_{\mathsf V_{n,m}}$ is a definable morphism. Moreover, we can prove that $h_{\phi}(\mathsf V_{n,m})$ and $h_{\phi}|_{\mathsf V_{n,m}}$ are small, and that the fiber of $h_{\phi}|_{\mathsf V_{n,m}}$ over $[(\varphi_1,\varphi_2,\varphi_3)]\in h_{\phi}(\mathsf V_{n,m})$ (of residue field $k'$) equals
$$\big\{\tau\in h_{\mathbb G_{m,k'(\!(t)\!)}} \mid -\ord_t\varphi_1\leq \ord_t\tau<\ord_t\varphi_2\big\}\cong \big\{\tau\in h_{\mathbb G_{m,k'(\!(t)\!)}} \mid 0\leq \ord_t\tau<m\big\}.$$
By \cite[Section 16.3]{CL}, $h_{\phi}|_{\mathsf V_{n,m}}$ gives rise to a $\mu_n$-equivariant semi-algebraic morphism of semi-algebraic sets in Denef-Loeser's framework $p:\mathscr{W}_{n,m}^{\infty} \to \widetilde{\mathscr{W}}_{n,m}^{\infty}$ with fiber over a point of residue field $k'$ isomorphic to $\{\tau\in \mathscr L(\mathbb A_{k'}^1) \mid 0\leq \ord_t\tau<m\}$. Finally, we can show that $p$ induces a $\mu_n$-equivariant constructible morphism of constructible sets $p_n:\mathscr{W}_{n,m}\to \widetilde{\mathscr{W}}_{n,m}$ with fiber 
$$\{\tau\in \mathscr L_n(\mathbb A_{k'}^1) \mid 0\leq \ord_t\tau<m\}\cong \mathbb A_{k'}^{n+1}\setminus \mathbb A_{k'}^{n+1-m},$$ 
as desired.



\end{document}